%% file: FacialBehav.tex
\documentclass[11pt]{article}
\usepackage{amssymb}
\usepackage{amsmath}
\usepackage{amsthm}
\newtheorem{theorem}[equation]{Theorem}
\newtheorem{corollary}[equation]{Corollary}
\newtheorem{definition}[equation]{Definition}
\newtheorem{proposition}[equation]{Proposition}

\newtheorem{remark}[equation]{Remark}
\newcommand\beq{\begin{equation}}

\newcommand\eeq{\end{equation}}
\newcommand\re{\mathrm {Re~}}

\newcommand\im{\mathrm {Im~}}
\newcommand\ii{\mathrm i}
\newcommand\nt{\stackrel{\rm  nt }{\to}}
\newcommand\al{\alpha}
\newcommand\de{\delta}
\newcommand\ga{\gamma}
\newcommand\ep{\varepsilon}
\newcommand\la{\lambda}
\renewcommand\ln{{\lambda_n}}

\newcommand{\ph}{\varphi}
\newcommand\si{\sigma}

\newcommand\D{\mathbb D}
\newcommand\T{\mathbb T}
\newcommand\C{\mathbb C}
\newcommand\M{\mathcal{M}}
\newcommand\R{\mathbb R}
\newcommand\e{\mathrm e}
\newcommand\Pick{\mathcal P}
\newcommand{\PP}{\mathcal{P}_2}
\newcommand\Schur{\mathcal{S}}
\renewcommand{\SS}{\mathcal{S}_2}

\newcommand\lan{\langle}
\newcommand\ran{\rangle}

\newcommand\half{{\tfrac{1}{2}}}

\newcommand\df{\stackrel{\rm def}{=}}
\newcommand\nn{\nonumber}
\newcommand\tb{\partial (\D^2)} 
\newcommand\tbR{\partial (\R^2)} 
\renewcommand\phi{\varphi}
\renewcommand\epsilon{\varepsilon}

\numberwithin{equation}{section}

\begin{document}
\title {Facial behaviour of analytic functions on the bidisk}
\author{Jim Agler, John E. McCarthy and N. J. Young}

\bibliographystyle{plain}
\date{\empty}
\maketitle
\begin{abstract}
We prove that if $\ph$ is an analytic function bounded by $1$ on the bidisk $\D^2$ and $\tau$ is a point  in a face of $\D^2$ at which $\ph$ satisfies Carath\'eodory's condition then both $\ph$ and the angular gradient $\nabla\ph$ exist and are constant on the face.  Moreover, the class of all $\ph$ with prescribed $\ph(\tau)$ and $\nabla\ph(\tau)$ can be parametrized in terms of a function in the two-variable Pick class.  As an application we solve an interpolation problem with nodes lying on faces of the bidisk.\footnote{Agler was partially supported by National Science Foundation Grant
DMS 0801259. M\raise.5ex\hbox{c}Carthy
was partially supported by National Science Foundation Grant DMS 0501079. Young
was supported by EPSRC Grant EP/G000018/1. \\
\hspace*{0.5cm} MSC Classification:  32A40 (primary), 32A07, 30J99, 32A30  (secondary).}
\end{abstract}
\input intro
\input notation
\input facial

\input cayley
\input julia
\input param
\input applic

\input FacialBehav_bbl
\noindent{Jim Agler,
Department of Mathematics, U.C. San Diego, La Jolla, CA 92093, USA}\\

\noindent{John E. M\raise.5ex\hbox{c}Carthy,  Department of Mathematics,  Washington University, St. Louis, MO 63130, USA  and Trinity College, Dublin, Ireland}\\

\noindent{N. J. Young,
 Department of Pure Mathematics, Leeds University, Leeds LS2 9JT, England.}
\end{document}

%% file: intro.tex
\section{Introduction} \label{intro}
We study functions in the Schur class $\SS$ of the bidisk $\D^2$, that is, functions analytic on $\D^2$ and bounded in absolute value by $1$, and in particular their behaviour on faces of   $\D^2$.  A {\em face} of $\D^2$ is a subset of the topological boundary $\tb$ of $\D^2$ having one of the forms $\{\tau^1\}\times\D$ or $\D\times\{\tau^2\}$, where $|\tau^1|=|\tau^2|=1$ and $\D$ denotes the open unit disk in the complex plane $\C$.   Of course a function $\ph\in\SS$ need not have values at all points of $\tb$, but it is known  \cite{jaf93,ab98} that $\ph$ has a nontangential limit at any point $\tau\in\tb$ for which the Carath\'eodory condition
\beq \label{cc}
\liminf_{\la\to\tau}\frac{1-|\ph(\la)|}{1-\|\la\|} < \infty
\eeq
holds, where 
\[
\| (\la^1,\la^2)\| = \max \{|\la^1|, |\la^2|\}.
\]
 We shall say that $\tau\in\tb$ is a {\em $B$-point for} $\ph$ if condition (\ref{cc}) holds.  We denote the nontangential limit (explained in detail below) of $\ph$ at a $B$-point $\tau$ by $\ph(\tau)$.

The purpose of this paper is to show that if $\ph\in\SS$ has a $B$-point $\tau$ lying on a face of $\D^2$ then strong consequences follow:  if $|\tau^1|=1$ then both $\ph$ and  the angular gradient $\nabla\ph$ exist and are constant on the face, and in fact
\[
\nabla\ph(\si) = \ph(\tau)\overline{\tau^1}\begin{pmatrix} \al \\ 0 \end{pmatrix}
\]
for all $\si\in\{\tau^1\}\times\D$, where $\al$ is the value of the $\liminf$ on the left hand side of inequality (\ref{cc}).
  Moreover, for any given   $\tau\in\T\times\D$,  the functions in $\ph\in\SS$ having a $B$-point at $\tau$ and prescribed values of $\ph(\tau)$ and $\nabla\ph(\tau)$ can be parametrized in terms of a function in the two-variable Pick class.

In Section \ref{notation} we establish notation and introduce the notion of angular gradient for functions in the two-variable Schur and Pick classes.  We also explain the notion of a {\em model} of a function in the Schur class and recall one property of models.  In Section \ref{facial} we prove the constancy result mentioned above (Theorem \ref{Bfaces}).  In Section \ref{cayley} we present the precise relationship between functions in the Schur and Pick classes given by the Cayley transform and thereby obtain an analog of Theorem \ref{Bfaces} for the Pick class.  In Section \ref{julia} we present Julia's reduction method for the Pick class, and in Section \ref{param} we derive two parametrization results, Theorems \ref{paramP2} and  \ref{paramS2}.   In Section \ref{applic} we apply our parametrization theorem to solve an interpolation problem.

%% file: notation.tex
\section{Notation and definitions} \label{notation}
We denote by $\Pi$ the upper half-plane $\{z\in\C: \im~z > 0\}$.  The one- and two-variable {\em Pick classes}, denoted by $\Pick$ and $\PP$, are the sets of analytic functions on $\Pi$ and $\Pi^2$ respectively with non-negative imaginary part.  The one- and two-variable {\em Schur classes} $\Schur$ and $\SS$ are the the sets of analytic functions on $\D$ and $\D^2$ respectively that are bounded by $1$ in absolute value.  We denote the closure of a set $S$ by $S^-$.

We shall need the notion of {\em nontangential approach}.  For a domain $\Omega\subset\C^n$, a set $S\subset \Omega$ is said to approach $\tau\in\partial\Omega$ nontangentially if $\tau\in S^-$ and
\[
 \frac {\|\la-\tau\|}{\mathrm{dist}(\la,\C^n\setminus \Omega)} \mbox{ is bounded for } \la\in S.
\]
The smallest $c \geq 1$ that bounds the above set is called the {\em aperture} of $S$.

In one variable the Julia-Carath\'eodory Theorem \cite{car29,car54,sar93} tells us that a function $\ph\in\Schur$ has an angular derivative at any $B$-point $\tau$.  More fully, if the Carath\'eodory condition holds (replace $\|\la\|$ by $|\la|$ in condition (\ref{cc})), then
\begin{enumerate}
\item[\rm (1)] $\ph(\la)$ tends to a limit $\ph(\tau)$ as $\la\to \tau$ in any set that approaches $\tau$ nontangentially;
\item[\rm (2)] the difference quotient $(\ph(\la)-\ph(\tau))/(\la - \tau)$ tends to a limit $\ph'(\tau)$ as $\la\to \tau$ in any set that approaches $\tau$ nontangentially;
\item[\rm(3)] $\ph'(\la)\to \ph'(\tau)$ as $\la\to\tau$ nontangentially;
\item[\rm(4)] the limit inferior in relation (\ref{cc}) is equal to the $\liminf$ along the radius $\la = r\tau, \ 0\leq r < 1$, and is in fact a limit as $r\to 1$.
\end{enumerate}
As we have mentioned,  (1) remains true for the polydisk \cite{jaf93,ab98}, and so $\ph\in\SS$ has a value $\ph(\tau)\in\T$ at any $B$-point $\tau$.  However, the analogue of (2) is false in two variables: in general $\ph\in\SS$ does not have an angular gradient at all its $B$-points  (see for example Remark \ref{counterex} below).  We are led to study those points $\tau$ at which $\ph$ does have an angular gradient.
\begin{definition} \label{4.4}
Let $\phi \in \mathcal{S}_2, \tau \in \partial(\D^2)$.\\
{\rm (1)} For $S \subseteq \mathbb{D}^2$, $\tau \in S^-$ we say that $\phi$ {\em  has a holomorphic differential at $\tau$ on $S$} if there exist $\omega, \eta^1, \eta^2 \in \mathbb{C}$ such that, for all $\lambda \in S$,
\begin{equation} \label{4.5}
\phi(\lambda) = \omega + \eta^1(\lambda^1-\tau^1) + \eta^2(\lambda^2 - \tau^2) + e(\lambda) 
\end{equation}
where 
\begin{equation} \label{4.6}
\lim_{\lambda\rightarrow \tau, \, \, \lambda \in S} \frac{e(\lambda)}{||\lambda-\tau||} = 0.
\end{equation} 
{\rm(2)} We say that $\tau$ is  a {\em $C$-point for}  $\phi$  if, for every set $S$ that approaches $\tau$ nontangentially, $\ph$ has a holomorphic differential on $S$ and $\omega$ in relation {\rm(\ref{4.5})} has modulus $1$. \\
{\rm(3)} If $\tau$ is a $C$-point for $\ph$ we define the {\em angular gradient $\nabla\ph (\tau)$} of $\ph$ at $\tau$ to be the vector $\begin{pmatrix} \eta^1 \\ \eta^2\end{pmatrix}$, where $\ph$ has holomorphic differential {\rm (\ref{4.5})} on some set that approaches $\tau$ nontangentially.
\end{definition}
It is clear that, when $\tau$ is  a  $C$-point for $\ph$, the quantities $\omega, \eta^1,\eta^2$ in equation (\ref{4.5}) are the same for every nontangential approach region $S$, and so the definition of $\nabla\ph (\tau)$ in (3) is unambiguous. 

An apparent drawback of  the above definition of $C$-point is that a condition must hold for {\em every} set $S$ that  approaches $\tau$ nontangentially.  However, we showed in \cite[Remark 8.12]{amy}  that the condition need only be checked for a single suitable set $S$.

Analogous definitions of $C$-points and angular gradients apply to functions in the Pick class $\PP$.

Our approach makes use of the notion of a model of a Schur-class function, as developed in \cite{amy}.
\begin{definition} \label{defmodel}
Let $\phi \in \mathcal{S}$. We say that $(\mathcal{M}, u)$ is a {\em model of} $\phi$ if $\mathcal{M} = \mathcal{M}^1 \oplus \mathcal{M}^2$ is an orthogonally decomposed separable Hilbert space and $u: \mathbb{D}^2 \rightarrow \mathcal{M}$ is an analytic map such that, for all $\lambda, \mu \in \mathbb{D}^2$,
\begin{equation}\label{modeleq}
1 - \overline{\phi(\mu)}\phi(\lambda) = (1 - \overline{\mu^1} \lambda^1) \langle  u_\lambda^1,u_\mu^1
\rangle  + (1 - \overline{\mu^2} \lambda^2)\langle u_\lambda^2,u_\mu^2 \rangle.  
\end{equation}
\end{definition}
In equation (\ref{modeleq}) we have written $u_\la$ for $u(\lambda)$, $u_\lambda^1 = P_{\mathcal{M}^1} u_\lambda$, and $u_\lambda^2 = P_{\mathcal{M}^2} u_\lambda$.   More generally, if $v \in \mathcal{M}$, we set $v^1 = P_{\mathcal{M}^1} v$ and $v^2 = P_{\mathcal{M}^2} v$. If $\lambda \in \mathbb{D}^2$, we may regard $\lambda$ as an operator on $\mathcal{M}$ by letting 
\[ 
\lambda v = \lambda^1 v^1 + \lambda^2 v^2
\] 
for $v \in \mathcal{M}$. 

Our methods depend on the fact that every function in $\SS$ has a model in the sense of Definition \ref{defmodel}, as was proved in \cite{ag90}.

If $(\M,u)$ is a model of $\ph\in\SS$ then we define the {\em nontangential cluster set} $X_\tau$ of the model at a $B$-point $\tau$ of $\ph$ to be the set of weak limits of weakly convergent sequences $u_{\la_n}$ over all sequences $\ln$ that converge nontangentially to $\tau$ in $\D^2$.    Here is a simple but powerful consequence of the model relation (\ref{modeleq}).
\begin{proposition}\label{getx}
Let $\tau\in\T\times\D$ be a $B$-point for $\ph\in\SS$ and let $(\M,u)$ be a model of $\ph$.  Then
\begin{enumerate}
\item[\rm (1)] $\emptyset \neq X_\tau \subset \M^1$;
\item[\rm (2)] there exists $\omega\in\T$ such that, for every $x=x^1 \in\ X_\tau$ and $\la\in\D^2$,
\beq\label{propx}
1-\overline{\omega} \ph(\la) = (1-\overline{\tau^1}\la)\lan u_\la^1, x^1 \ran.
\eeq
\end{enumerate}
\end{proposition}
A detailed proof is given in \cite[Proposition 4.2]{amy}.  Nonemptiness of $X_\tau$ follows from the fact that $u_{\la_n}$ is bounded for any sequence $(\ln)$ in $\D^2$ that converges to 
$\tau$ nontangentially, while the relation (\ref{propx}) is derived by substituting $\mu=\mu_n$ in equation (\ref{modeleq}) for a suitable sequence $(\mu_n)$ converging to $\tau$ and then passing to a limit.

It is shown in \cite[Corollary 8.11]{amy} that $\tau$ is a $C$-point for $\ph$ if and only if $X_\tau$ is a singleton set; we denote the unique member of this set by $u_\tau$.

%% file: facial.tex
\section{Facial $B$-points} \label{facial}
 A function $\ph\in\SS$ can have a $B$-point in $\partial (\D^2) \setminus \T^2$, say at $\tau\in\T\times \D$; we shall call such a $\tau$ a {\em facial $B$-point} for $\ph$.  Facial $B$-points can arise in a trivial way, when $\ph$ is independent of $\la^2$, but can also occur non-trivially.  Consider for example the rational inner function
\beq\label{ratex}
\ph(\la) = \frac{1+\la^1+\la^2-3\la^1\la^2}{3-\la^1-\la^2-\la^1\la^2}.
\eeq
$\ph$ is analytic and equal to $1$ at every point of the face $\{1\}\times \D$ of the bidisk, and so every point of this face is a $C$-point for $\ph$.  Likewise $\ph=1$ at every point of the face $\D\times\{1\}$.
The example illustrates a general phenomenon.  We shall denote by $\Delta$ the closed unit disk $\{z\in\C: |z|\leq 1\}$.

\begin{theorem} \label{Bfaces}
Let $\tau\in\T\times\D$ be a $B$-point for $\ph\in\SS$.   Then 
\begin{enumerate}
\item [\rm (1)] every point of $\{\tau^1\}\times \Delta$ is a $B$-point for $\ph$;
\item[\rm (2)]  every point of $\{\tau^1\}\times \D$ is a $C$-point for $\ph$;
\item [\rm (3)]   $\ph$  is constant on $\{\tau^1\}\times\Delta$; 
\item [\rm (4)]    $\nabla\ph$ is constant on $\{\tau^1\}\times \D$, with value 
\beq \label{valnablaph}
\nabla\ph(\si) = \ph(\tau)\overline{\tau^1}\al \begin{pmatrix} 1 \\ 0 \end{pmatrix} \quad \mbox{ for every } \si\in \{\tau^1\}\times \D
\eeq
where
\[
\al =  \liminf_{\la\to\tau} \frac{1-|\ph(\la)|}{1-\|\la\|}.
\]
\end{enumerate}
\end{theorem}
\begin{proof}
Let $(\M,u)$ be a model of $\ph$.  
By Proposition  \ref{getx}, since $|\tau^2| < 1$, there exists $x\in \M$ and $\omega\in\T$ such that $x^2=0$ and
\beq\label{formphi}
1-\bar\omega\ph(\la) = (1-\overline{\tau^1}\la^1) \lan u_\la^1,x^1 \ran
\eeq
and so
\beq \label{j2a}
1-|\ph(\la)| \leq |1-\overline\omega\ph(\la)| \leq |1-\overline{\tau^1}\la^1| \ | \lan u_\la^1,x^1 \ran |
\eeq
for all $\la\in\D^2$.\\
(1)  Consider any point $\si=(\tau^1,\si^2) \in \{\tau^1\}\times\D$.   Let
\[
S_\si = \{((1-t)\tau^1, \si^2): 0<t\leq 1-|\si^2|\}.
\]
$S_\si$  approaches $\si$ nontangentially with aperture $1$, and we have
\beq\label{Ssieqal}
|1-\overline{\tau^1}\la^1| = 1-\|\la\|
\eeq
for $\la\in S_\si$.
 By \cite[Remark 5.6]{amy}, for $\la\in S_\si$,
\beq \label{boundu1}
||u_\la|| \leq 2||x^1||,
\eeq
which, together with (\ref{j2a}) and (\ref{Ssieqal}), implies that
\beq\label{Cquo}
\frac{1-|\ph(\la)|}{1-||\la||} = \frac{1-|\ph(\la)|}{ |1-\overline{\tau^1}\la^1|} \leq 2 ||x^1||^2
\eeq
for all $\la\in S_\si$.  Hence $\si$ is a $B$-point for $\ph$.

Now consider any point $\rho\in \{\tau^1\}\times\T$.  For any $r\in(0,1)$ observe that
\[
S_{(\tau^1,(1-r)\rho^2)} = \{((1-t)\tau^1,(1-r)\rho^2): 0<t\leq r\},
\]
and so in particular $(1-t)\rho\in S_{(\tau^1, (1-t)\rho^2)}$.  The bound (\ref{boundu1}) is valid for all $\la\in\cup_{\si\in\{\tau^1\}\times\D} S_\si$ and hence for $\la=(1-t)\rho$.  We therefore have the bound (\ref{Cquo}) for $\la=(1-t)\rho$ which tends to $\rho$ as $t\to 0+$.   Hence $\rho$ is a $B$-point for $\ph$.\\

\noindent (2)   According to \cite[Lemma 8.10]{amy}, the nontangential cluster set of $u$ at any facial $B$-point $\si$ comprises a single vector in $\M$, and so every facial $B$-point is a $C$-point for $\ph$.\\

\noindent (3) Equation (\ref{formphi}) can be re-written
\[
\ph(\la) = \omega +\omega\overline{\tau^1}(\la^1-\tau^1) \lan u_\la^1, x^1 \ran.
\]
Consider $\si\in\{\tau^1\}\times \D$.  Substitute $\ln$ for $\la$ and let $\ln \to \si$ in $S_\si$; by inequality (\ref{boundu1}), the second term on the right hand side tends to zero, and we find that 
\[
\ph(\si) \df \lim_{\la \nt \si} \ph(\la) = \omega.
\]
For $\rho\in\{\tau^1\}\times\T$ we may substitute $\la=\ln=(1-t_n)\rho$ where $t_n \to 0+$, and similar reasoning applies.
Thus $\ph$ is constant on $\{\tau^1\}\times \Delta$.\\

\noindent (4) We prove that $\nabla\ph$ is also constant.  For any $\si=(\tau^1, \zeta), \ \zeta\in\D$, let
\[
\nabla\ph(\si)= \begin{pmatrix} \eta^1(\zeta) \\ \eta^2(\zeta) \end{pmatrix}.
\]
By \cite[Corollary 8.13]{amy} and the fact that $u_\si^2=0$ we have
\beq \label{poscomp}
  \begin{pmatrix} \eta^1(\zeta) \\ \eta^2(\zeta) \end{pmatrix} =\nabla\ph(\si) =\ph(\si)\begin{pmatrix}\overline{\si^1}|| u^1_\si ||^2 \\
   \overline{\si^2}|| u^2_\si ||^2 \end{pmatrix} = \ph(\tau)\overline{\tau^1}\begin{pmatrix} \|u_\si^1\|^2 \\ 0\end{pmatrix},
\eeq
and hence 
\beq \label{poseta}
\overline{\ph(\tau)}\tau^1 \eta^1(\zeta) \geq 0 \mbox{ for all } \zeta\in\D.
\eeq

 The points $(\tau^1-t\tau^1, \zeta), \ 0< t\leq 1,$ approach $\si$ nontangentially, and so we have
\[
\eta^1(\zeta) =  \lim_{t\to 0+}F_t(\zeta)
\]
where
\[
F_t(\zeta) = -\overline{\tau^1} \frac{\ph(\tau^1-t\tau^1,\zeta)-\ph(\tau)} {t}.
\]
We claim that the functions $F_t, \ 0 < t \leq 1,$ are uniformly bounded on $\D$.  Indeed, by equation (\ref{formphi}),
\[
\ph(\tau^1-t\tau^1, \zeta)-\ph(\tau) = -\ph(\tau)t\lan u_{(\tau^1-t\tau^1, \zeta)}, x^1 \ran
\]
and therefore
\[
|F_t(\zeta)| \leq \| u_{(\tau^1-t\tau^1, \zeta)}\| \ \| x^1\|.
\]
The set $\{ (\tau^1-t\tau^1, \zeta) : 0 < t \leq 1\}$ approaches $\si$ with aperture $1$, and so, by \cite[Remark 5.6]{amy},
\beq \label{boundu}
\|u_{(\tau^1-t\tau^1, \zeta)}\| \leq 2 \|x^1\|
\eeq
and consequently $|F_t(\zeta)| \leq 2\|x^1\|^2$ for $0<t\leq 1, \ \zeta\in\D$.  Since the $F_t$ are analytic in $\D$, so is their pointwise limit $\eta^1$.  In view of the positivity relation (\ref{poseta}), it follows that $\eta^1$ is constant on $\D$.
By equation (\ref{poscomp}) and \cite[Theorem 5.10]{amy} we therefore have
\beq \label{constnormu}
\|u_\si^1\|^2=\|u_\tau^1\|^2=\|u_\tau\|^2 = \al,
\eeq
and so the constant value of $\nabla\ph(\si)$ is given by equation (\ref{valnablaph}).
\end{proof} 

\begin{corollary} \label{2faces}
If a function $\ph\in \SS$ has a $B$-point in $\T\times\D$ and another in $\D\times\T$ then $\ph$ takes the same value at the two $B$-points.
\end{corollary}
For if $\tau\in\T\times\D$ and $\si\in\D\times\T$ are $B$-points then the constant value of $\ph$ on both of the closed faces must equal $\ph(\tau^1,\si^2)$.

\begin{remark}  \label{counterex} {\rm
It is not the case that if $\ph$ has a facial $B$-point $\tau\in\T\times \D$ then every point of the {\em closed} face $\{\tau^1\}\times\Delta$ is necessarily a $C$-point for $\ph$.
A counterexample is furnished by the rational inner function
\[
\psi(\la) = \frac{2\la^1\la^2-\la^1-\la^2}{2-\la^1-\la^2}.
\]
The face $\{1\}\times\D$ consists of $B$-points for $\psi$, but the point $(1,1)$ is not a $C$-point for $\psi$, and so $\nabla\psi(1,1)$ is undefined.  This example is analysed in \cite[Section 6]{amy}.}
\end{remark} 
A modification of the proof of Theorem \ref{Bfaces} yields a slightly stronger result.
\begin{proposition} \label{uconst}
Let $\tau\in\T\times\D$ be a $B$-point for $\ph\in\SS$.  For every model $(\M,u)$ of $\ph$,   for all $\si\in\{\tau^1\}\times\D$ we have  $u_\si=u_\tau$.
\end{proposition}
\begin{proof}
Choose $x\in\M$ as in the preceding proof.   For $\zeta\in\D$ let $\si=(\tau^1,\zeta)$.
We claim that $u_{(\tau^1-t\tau^1,\zeta)} \to u_\si$ as $t\to 0+$.  Indeed,  $u_{(\tau^1-t\tau^1,\zeta)}= u_{\si - t\de}$ where $\de=(\tau^1,0)$, and  by \cite[Theorem 7.1]{amy}, $u_{\si - t\de}$ tends to a limit in the cluster set $X_\si$ of $u$ at $\si$ as $t\to 0+$.  Since $\si$ is a $C$-point for $\ph$, \ $X_\si$ comprises the unique vector $u_\si$, which proves the claim.

Consider the functions $\zeta \mapsto u_{(\tau^1-t\tau^1,\zeta)}:\D\to\M, \ 0 < t \leq 1$.  It follows from inequality (\ref{boundu}) that these analytic functions are uniformly bounded on $\D$.  Since they are also pointwise convergent on $\D$ to the function $\zeta\mapsto u_{(\tau^1,\zeta)}= u_\si$, the latter map is analytic on $\D$.   By equation (\ref{constnormu}), $\|u_\si\|$ is constant on $\D$.  However, an analytic  Hilbert-space-valued function $f$ such that $\| f(.)\|$ is constant is itself a constant function.  Thus $u_\si$ is constant on $\{\tau^1\}\times\D$.
\end{proof}

%% file: cayley.tex
\section{The Cayley transform for $\SS$ and $\PP$} \label{cayley}
The Cayley transform 
\[
C: \D \to \Pi \mbox{ and } \C\setminus \{1\} \to\C\setminus \{-\ii\}: \la \mapsto \ii\frac{1+\la}{1-\la}
\]
enables us to pass back and forth between $\Schur$ and $\Pick$.  The transform $C$ also acts (co-ordinatewise) from $\D^2$ to $\Pi^2$.   The relationship between properties of a function $\ph\in \SS$ and those of the corresponding function $h\in\PP$ is straightforward, but not totally trivial, and so we summarize it here.

The variables $\la\in\D, \ z\in\Pi$ will be related by $z=C(\la)$, so that
\beq \label{laz}
z=\ii\frac{1+\la}{1-\la}, \qquad \la = \frac{z-\ii}{z+\ii}.
\eeq
Similarly, variables $\la\in\D^2$ and $z\in\Pi^2$ are related co-ordinatewise ($z^j=C(\la^j)$ etc).

Let $h\in\PP$ correspond to  $\ph\in\SS$, \ $\ph$ not identically equal to $1$,  under the Cayley transform, that is,
\beq \label{defh}
h(z) = \ii \frac{1+\ph(\la)}{1-\ph(\la)},  \qquad \ph(\la) = \frac{h(z)-\ii}{h(z)+\ii},
\eeq
where $\la, z$ are related by equations (\ref{laz}).  We consider a $B$-point $\tau\in\tb$ for $\ph$, and we let  $x\in \tbR$ correspond to $\tau$:
\beq \label{defx}
x^j =\ii\frac{1+\tau^j}{1-\tau^j}, \qquad \tau^j=\frac{x^j-\ii}{x^j+\ii}, \qquad \mbox{ for }  j=1,2.
\eeq
Let $\ph(\tau)=\omega$, so that  $|\omega|=1$.
Then $h(x)= \xi$ where $\xi\in\R\cup \{\infty\}$ is the transform of $\omega\in\T$:
\[
\xi=\ii\frac{1+\omega}{1-\omega}, \qquad \omega =\frac{\xi-\ii}{\xi+\ii}.
\]
We ask: what conditions on $h, x$ correspond to $\tau$ being a $B$-point or a $C$-point of $\ph$, and what is the relation between $\nabla\ph(\tau)$ and $\nabla h(x)$ in the case that $\tau$ is a $C$-point of $\ph$?

Let us assume that $\tau^1, \tau^2$ and $\omega$ are all different from $1$ (else we may compose with rotations of the circle).  It is then the case that $x^1, x^2$ and $\xi$ are all real numbers (not $\infty$).  By a straightforward calculation,
\[
\frac{1-|\ph(\la)|^2}{1-\|\la\|^2} = \frac{\im h(z)}{|h(z)+\ii|^2} \max_{j}\frac{|z^j +\ii|^2}{\im z^j}.
\]
Hence 
\beq \label{Climinf}
\liminf_{\la\to\tau}  \frac{1-|\ph(\la)|^2}{1-\|\la\|^2} = \liminf_{ z\to x} \frac{1}{|\xi+\ii|^2} \max_j \frac{|x^j + \ii|^2 \im h(z)}{\im z^j}.
\eeq
It follows that $\tau$ is a $B$-point for $\ph$ if and only if
\beq\label{hBpt}
\liminf_{z\to x} \frac{\im h(z)}{\im z^j} < \infty \mbox{ for } j=1,2.
\eeq
We shall say that $x\in\tbR$ is a $B$-point for $h\in\PP$ whenever the relation (\ref{hBpt}) holds.

Note that, for $\tau\in\T\times\D, \ x\in\R\times\Pi$, the limits inferior in equation (\ref{Climinf}) occur for $j=1$.  Hence, if
\beq \label{defal}
\al \df  \liminf_{\la\to\tau}  \frac{1-|\ph(\la)|}{1-\|\la\|}  = \liminf_{\la\to\tau}  \frac{1-|\ph(\la)|}{1-|\la^1|}
\eeq
and
\beq\label{defbeta}
\beta \df \liminf_{z\to x} \max_j \frac{\im h(z)}{\im z^j} = \liminf_{z\to x} \frac{\im h(z)}{\im z^1}
\eeq
then equation (\ref{Climinf}) yields
\beq\label{albeta}
\al = \left| \frac{x^1+\ii}{\xi+\ii}\right|^2\beta.
\eeq

Now consider a $C$-point $\tau$ for $\ph\in\SS$.   There exist $\omega\in\T$ and $\eta^1,\eta^2\in\C$ such that equation (\ref{4.5}) holds on any set $S$ that approaches $\tau$ nontangentially.  Nontangential approach is preserved by the Cayley transform, and
\[
\|\la-\tau\| \to 0 \Leftrightarrow \|z-x\| \to 0.
\]
Let us rewrite equation (\ref{4.5}) in terms of the variables $z,\ x$ and $h$.  We have
\begin{align} \label{hdiff}
h(z) &= \ii\frac{1+\ph(\la)}{1-\ph(\la)} \nn \\
    &= \ii \frac{1+\omega+ \eta\cdot(\la-\tau) +o(\|\la-\tau\|)}{1-\omega- \eta\cdot(\la-\tau) +o(\|\la-\tau\|)}  \nn \\
    &= \xi \frac{1+(1+\omega)^{-1} \eta\cdot(\la-\tau) +o(\|\la-\tau\|)}{1-(1-\omega)^{-1} \eta\cdot(\la-\tau) +o(\|\la-\tau\|) } \nn \\
    &= \xi\left\{ 1+\left(\frac{1}{1+\omega}+\frac{1}{1-\omega}\right)\eta\cdot(\la-\tau)+ o(\|\la-\tau\|) \right\} \nn \\
    &= \xi + \frac{2\xi}{1-\omega^2} \eta\cdot(\la-\tau) +o(\|\la-\tau\|)
\end{align}
as $\la\to \tau$ in $S$.  Now
\[
\la^j-\tau^j = \frac{2\ii(z^j-x^j)}{(z^j+\ii)(x^j+\ii)} = \frac{2\ii(z^j-x^j)}{(x^j+\ii)^2} + o{(\|z-x\|)}
\]
as $z\to x$.  Hence
\[
h(z) = \xi + \sum_{j=1,2} \left(\frac{1-\tau^j}{1-\omega}\right)^2\eta^j(z^j-x^j) + o(\|z-x\|)
\]
as $z\to x$ in any nontangential approach region.  That is, $x$ is a $C$-point for $h$, and
\beq\label{gotnablah}
\nabla h(x) = \frac{1}{(1-\omega)^2}\begin{pmatrix} (1-\tau^1)^2\eta^1\\ (1-\tau^2)^2\eta^2 \end{pmatrix}.
\eeq
In the case that $\tau$ is a facial $B$-point we have the following conclusion.
\begin{proposition} \label{valnablah}
Let $\tau\in\T\times\D$ be a $B$-point for $\ph\in\SS$ and suppose that $\tau^1\neq1, \ \ph(\tau)\neq 1$.  Then $x\in \R\times\Pi$ given by equations {\rm (\ref{defx})} is a $C$-point for $h\in\PP$ given by equations {\rm (\ref{defh})}, and
\[
\nabla h(x) = \left|\frac{1-\tau^1}{1-\omega}\right|^2  \begin{pmatrix} \al\\ 0 \end{pmatrix}
\]
where $\al$ is defined by equation  {\rm (\ref{defal})}.
\end{proposition}
For on combining equations (\ref{valnablaph}) and (\ref{gotnablah}) we find
\[
\nabla h(x) = \frac{(1-\tau^1)^2}{(1-\omega)^2} \omega \overline{\tau^1} \al \begin{pmatrix} 1 \\ 0 \end{pmatrix} = \left|\frac{1-\tau^1}{1-\omega}\right|^2 \al \begin{pmatrix} 1\\ 0 \end{pmatrix}.
\]
There is of course a converse to this result, which we shall not spell out.

We can derive a version of Theorem \ref{Bfaces} for the Pick class.  
\begin{theorem} \label{BfacesPick}
Let $x\in\R\times\Pi$ be a $B$-point for $h\in\PP$.  Then 
\begin{enumerate}
\item[\rm (1)] every point of $\{x^1\}\times \Pi^-$ is a $B$-point for $h$;
\item[\rm (2)] every point of $\{x^1\}\times \Pi$ is a $C$-point for $h$;
\item[\rm (3)] $h$  is constant on $\{x^1\}\times\Pi^-$;
\item[\rm(4)] $\nabla h$ is constant on $\{x^1\}\times\Pi$, with value
\beq \label{nablah}
\nabla h(y) =  \begin{pmatrix} \beta \\ 0 \end{pmatrix} \mbox{ for every } y\in \{x^1\}\times \Pi
\eeq
where
\[
\beta =  \liminf_{z\to x} \frac{\im h(z)}{\im z^1}.
\]
\end{enumerate}
\end{theorem}
\begin{proof}
(1)-(3) are immediate.  To check the value in (4)
apply Theorem \ref{Bfaces} to the function $\ph\in\SS$ defined by equation (\ref{defh}).
By Proposition \ref{valnablah} and equation (\ref{albeta}) we have
\[
\nabla h(y) = \nabla h(x) = \left|\frac{1-\tau^1}{1-\omega}\right|^2  \al\begin{pmatrix} 1\\ 0 \end{pmatrix} = \left|\frac{1-\tau^1}{1-\omega}\right|^2 \left |\frac{x^1+\ii}{\xi+\ii}\right|^2 \beta \begin{pmatrix} 1\\ 0 \end{pmatrix}.
\]
Since 
\[
(x^1+\ii)(1-\tau^1) = 2\ii = (\xi +\ii)(1-\omega),
\]
equation (\ref{nablah}) follows.
\end{proof}

%% file: julia.tex
\section{Julia reduction in $\Pick$} \label{julia}
G. Julia,  in the course of proving his celebrated Lemma in  \cite{ju20}, introduced a form of reduction for functions in the Pick class and showed that reduction preserves the Pick class.  His reduction process is analogous to the better-known Schur reduction for functions in the Schur class, but is associated with {\em boundary} points, that is, points on the real axis.  Julia's reduction was subsequently used extensively by R. Nevanlinna, e.g. in \cite{Nev1,Nev2}.  In the next section we shall apply it to functions in $\PP$ having a facial $B$-point.

Recall that, by the (one-variable) Julia-Carath\'eodory Theorem, if $x\in\R$ is a $B$-point for $f\in\Pick$ then $f$ has an angular derivative $f'(x)$ at $x$.   Furthermore, if $f$ is non-constant,
\[
f'(x) = \liminf_{z\to x}\frac{\im f(z)}{\im z}  = \lim_{y\to 0+} \frac{\im f(x+\ii y)}{y}> 0.
\]
\begin{definition} \label{defreduce} \rm
 (1)  For any non-constant function $f\in\Pick$ and any $x\in\R$ such that $x$ is a $B$-point for $f$ we define the {\em reduction of $f$ at $x$} to be the function $g$ on $\Pi$ given by the equation
\beq \label{reducef}
g(z) = -\frac{1}{f(z)-f(x)} + \frac{1}{f'(x)(z-x)}.
\eeq

 (2)  For any $g\in\Pick$, any $x\in\R$  and any $a_0\in\R, a_1 > 0$, we define the {\em augmentation of $g$ at $x$ by $a_0, a_1$} to be the function $f$ on $\Pi$ given by
\beq \label{augmentg}
 \frac{1}{f(z)-a_0} =  \frac{1}{a_1(z-x)} -g(z).
\eeq
\end{definition}
Note that in (1), since $f(x)~$ is real and $f$ is non-constant, the denominator $f(z)-f(x)$ is non-zero by the maximum principle.
Note also that $f$ defined by equation (\ref{augmentg}) is necessarily non-constant, for otherwise
\[
\im g(z) = \mathrm{const} + \frac{1}{a_1}\im \frac{1}{z-x},
\]
and the last term can be an arbitrarily large negative number for $z\in\Pi$, contrary to the choice of $g\in\Pick$.

Reduction and augmentation are of course inverse operations.

The following is a refinement of Julia's result.
\begin{theorem} \label{propfg}
Let $x \in\R$.
\begin{enumerate}
\item[\rm(1)] If  $x$ is a $B$-point for a non-constant function $f\in\Pick$ then the reduction of $f$ at $x$ also belongs to $\Pick$.
\item[\rm(2)] If $g\in\Pick$ and $a_0\in\R,\, a_1>0$ then the augmentation $f$ of $g$ at $x$ by $a_0,\, a_1$  belongs to $\Pick$, has a $B$-point at $x$  and satisfies $f(x)=a_0, \ f'(x) \leq a_1$.
Moreover
\[
 f'(x) = a_1  \quad\mbox{ if and only if }\quad \lim_{y\to 0+} yg(x+\ii y) =0.
\]
\end{enumerate}
\end{theorem}
\begin{proof}
  Julia proved (1) in the case that $f$ is regular at $x$ (and observed in a footnote that it is true slightly more generally).  The following proof of the general case is essentially due to Nevanlinna \cite[pp. 6-9]{Nev2}.

By the Julia-Carath\'eodory Theorem $f$ and $f'$ have nontangential limits $a_0,\ a_1$  respectively at $x$, and $a_0\in\R, \ a_1 > 0$.  Moreover, $x$ is a $C$-point for $f$, so that
\beq \label{Cptx}
f(z) = a_0+a_1 (z-x) + R(z) \quad \mbox{ for } z \in\Pi
\eeq
where
\[
\lim_{z \stackrel{\rm nt}{\to} 0} \frac{R(z)}{z-x} = 0.
\]

Let $g$ be the reduction of $f$ at $x$, so that 
\[
g(z) = -\frac{1}{f(z)-a_0} +\frac{1}{a_1(z-x)}.
\]
We have 
\[
\im \left( g(z) -\frac{1}{a_1(z-x)}\right) = \im \left(-\frac{1}{f(z)-a_0}\right) \geq 0.
\]
Let $w = (-1)/(z-x)$, so that $w\in\Pi$ if and only if $z\in\Pi$.  The last inequality can be written
\beq\label{imgz}
-\im g(z) \leq \frac{\im w}{a_1}.
\eeq

Let $\ep\in (0, \half\pi)$ and let
\[
U = \{r\mathrm{e}^{\ii\theta} \in\Pi: r>0, 0<\theta \leq \ep \mbox{ or } \pi-\ep \leq \theta < \pi\}, \quad V= \Pi \setminus U.
\]
Notice that $z\in U+x$ if and only if $w\in U$.  If $w=r\mathrm{e}^{\ii\theta}\in U$ then
\[
\im w = r\sin \theta\leq r\theta\leq \ep |w|
\]
and so, in view of the relation (\ref{imgz}),
\beq\label{imgonU}
-\im g(z) \leq \frac{\ep |w|}{a_1}  \quad \mbox{ for all } z\in U.
\eeq
We claim that the same inequality holds for  $w\in V$ of sufficiently large norm.  For
\[
g(z) = \frac{f(z)-a_0-a_1(z-x)}{a_1(z-x)(f(z)-a_0)} = \frac{R(z)}{a_1(z-x)(a_1(z-x)+R(z))}.
\]
Hence 
\[
(z-x)g(z) = \frac{R(z)/(z-x)}{a_1(a_1+R(z)/(z-x))} \to 0 \quad \mbox{ as } z \stackrel{\rm nt}{\to} x.
\]
Thus, for some $r_0>0$,
\[
-\im g(z) \leq |g(z)| < \frac{ \ep |w|}{a_1} \quad \mbox{ for all } w\in V, |w| > r_0.
\] 
The inequality (\ref{imgonU}) therefore holds for all $w\in\Pi$ such that $|w| > r_0$.

Define an analytic function $F$ on $\Pi$ by
\[
F(w) = \e^{\ii g(z)} = \e^{\ii g(x-1/w)}.
\]
We have, for any $w \in \Pi$,
\[
|F(w)| = \e^{\re \ii g(z)} = \e^{-\im g(z)},
\]
and hence, by (\ref{imgonU}), whenever $|w| > r_0$,
\[
|F(w)|  \leq  \e^\frac{\ep |w|}{a_1},
\]
that is, $F$ has at worst exponential growth on $\Pi$.  We may therefore apply the Phragm\'en-Lindel\"of Theorem to $F$ on the half-plane $\{w: \im w \geq \de\}$ for any $\de > 0$.  By the inequality (\ref{imgz}), if $\im w =\de$ then $-\im g(z) \leq \de/a_1$ and therefore
\[
|F(w)| = \e^{-\im g(z)} \leq \e^{\de/a_1}.
\]
It follows by Phragm\'en-Lindel\"of (e.g. \cite[p. 218]{BakNew}) that  $|F| \leq \e^{\de/a_1}$ on $\Pi+\ii\de$.  On letting $\de$ tend to zero we deduce that $|F| \leq 1$ on $\Pi$, and hence that $\im g \geq 0$ on $\Pi$.
Thus $g\in\Pick$.

To prove (2) consider any $g\in\Pick, \ a_0\in\R, \ a_1>0$ and let $f$ be the corresponding augmentation of $g$, so that equation (\ref{augmentg}) holds.
For any $y>0$,
\begin{align} \label{formratio}
\frac{\im f(x+\ii y)}{y} &= \frac{1}{y} \im\left( a_0 + \frac{1}{\frac{1}{a_1\ii y} - g(x+\ii y)}\right)  \nn \\
   &= a_1 \re \frac{1}{1-\ii a_1 y g(x+\ii y)}.
\end{align}
For any $c\in(0,\infty)$ and $z\in\C$ we have
\begin{align}\label{elobs}
\re \frac{1}{z} \leq c  \qquad &\Leftrightarrow   \qquad \left| z - \frac{1}{2c} \right| \geq \frac{1}{2c} \nn \\
   &\Leftrightarrow \qquad  z \notin {D}\left( \frac{1}{2c}, \frac{1}{2c}\right),
\end{align}
where ${D}(w,r)$  denotes the open disk of centre $w$ and radius $r$.

Now for any $y>0$ we have
\[
\re( 1 - \ii a_1 y g(x+\ii y)) = 1+ a_1y\im g(x+\ii y) \geq 1,
\]
and hence
\[
1 - \ii a_1 y g(x+\ii y) \notin {D}( \tfrac 12, \tfrac 12).
\]
It follows that
\[
\re \frac{1}{1-\ii a_1 y g(x+\ii y)} \leq 1
\]
and hence, by equation (\ref{formratio}), that
\[
\frac{\im f(x+\ii y)}{y} \leq a_1 \qquad \mbox{ for all } y>0.
\]
Thus $x$ is a $B$-point for $f$.

For $y>0$ we have
\[
\frac{1}{f(x+\ii y)- a_0} = -\frac{1}{a_1} \left( \frac{\ii}{y} + a_1 g(x+\ii y) \right).
\]
Now
\[
\im \left( \frac{\ii}{y} + a_1 g(x+\ii y) \right) = \frac{1}{y} + a_1\im g(x+\ii y) \geq \frac{1}{y},
\]
and hence
\[
\frac{1}{|f(x+\ii y) - a_0|} \to \infty \quad \mbox{ as } y \to 0+.
\]
Thus $f(x+\ii y) \to a_0$ as $y\to 0+$, which is to say that $f(x) = a_0$.

Again by equation (\ref{augmentg}) we have, for $y> 0$,
\[
\frac{\ii y}{f(x+\ii y)-a_0} = \frac{1}{a_1} - \ii yg(x+\ii y).
\]
The left hand side tends to the limit $1/f'(x)$ as $y\to 0+$ and hence the limit of $\ii yg(x+\ii y)$ exists and is real.  Now
\[
\re \ii y g(x+\ii y) =  -y\im g(x+\ii y) \leq 0,
\]
and hence 
\beq \label{fprime}
\frac{1}{f'(x)} =\frac{1}{a_1} + \lim_{y \to 0+} y\im g(x+\ii y) \geq \frac{1}{a_1}.
\eeq
It follows that $f'(x) \leq a_1$ and that $f'(x)=a_1$ if and only if 
\beq \label{nopole}
\lim_{y\to 0+} yg(x+\ii y) = 0.
\eeq
\end{proof}
 Examples of functions $g(z)$ in $\Pick$ for which formula (\ref{nopole}) holds are  $(z-x)^\al$ for $0\leq \al \leq 1$, $ -(z-x)^{-\al}$ for $0<\al<1$ and $\log (z-x)$.
Functions in $\Pick$ for which it does {\em not} hold are $-1/(z-x)$ and $-\cot(z-x)$.  Roughly speaking, condition (\ref{nopole}) rules out poles of $g$ at $x$. 

%% file: param.tex
\section{Parametrization theorems}\label{param}
The first theorem describes the functions in $\PP$ with a prescribed facial $B$-point.
\begin{theorem} \label{paramP2}
Let $x\in\R\times\Pi, \ \xi\in\R$ and $\beta>0$.  The functions $h\in\PP$ such that $x$ is a $B$-point for $h$,  $h(x)=\xi$ and $\nabla h(x) = (\beta, 0)^T$ are precisely the functions of the form
\beq\label{paramh}
h(z) = \xi + \frac{1}{\frac{1}{\beta(z^1-x^1)} - g(z)}
\eeq
for some function $g\in\PP$ such that
\beq\label{limyg}
\lim_{y\to 0+} yg(x^1+\ii y, x^2) = 0.
\eeq
\end{theorem}
\begin{proof}
Consider a function $h$ of the form (\ref{paramh}) for some $g$ as described.   For fixed $z^2\in\Pi$ it is clear that $g(.,z^2)$ is in the one-variable Pick class $\Pick$ and $h(.,z^2)$ is the augmentation of $g(.,z^2)$ at $x^1$ by $\xi,\ \beta$.  By Theorem \ref{propfg}(2), $h(.,z^2)\in\Pick, \ h(.,z^2)$ has a $B$-point at $x^1$ and
\[
\lim_{y\to 0+}  h(x^1+\ii y, z^2) = \xi.
\]
Since $(x^1+\ii y, z^2)$ tends to the $B$-point $(x^1,z^2)$ nontangentially as $y\to 0+$ we have $h(x^1,z^2) = \xi$.  Again by Theorem \ref{propfg}, the angular derivative of $h(.,z^2)$ at $x^1$ is at most $\beta $, and because of equation (\ref{limyg}), when $z^2=x^2$, this angular derivative is exactly $\beta$ at $x^1$.  It follows that $\nabla h(x) = (\beta, 0)^T$.

Conversely, suppose that $x$ is a $B$-point for $h\in\PP$, that $h(x)=\xi$ and that $\nabla h(x) = (\beta, 0)^T$.  By Theorem \ref{BfacesPick}, all points $(x^1,w)$ with $w\in\Pi$ are $C$-points for $h$ and we have $h(x^1,w)=\xi, \ \nabla h(x^1,w) = (\beta, 0)^T$.  For any $w\in\Pi$ it follows that $h(.,w)$ is a non-constant function in the Pick class, $x^1$ is a $C$-point for $h(.,w)$ and the value and angular derivative of $h(.,w)$ at $x^1$ are $\xi$ and $\beta$ respectively.   We may therefore reduce $h(.,w)$ at $x^1$ to obtain $g(.,w)\in\Pick$ given by
\[
g(z^1,w) = -\frac{1}{h(z^1,w)-\xi} + \frac{1}{\beta(z^1 - x^1)}.
\]
Clearly $g$ is analytic on $\Pi^2$, and so $g\in\PP$.  Now $h$ is related to $g$ by equation (\ref{paramh}).  Furthermore,
\[
\frac{1}{\beta}=\frac{1}{D_1 h(x^1,w)} = \frac{1}{\beta} + \lim_{y \to 0+} yg(x^1+\ii y,w)
\]
where $D_1$ here denotes the angular derivative in the first variable.  Hence
\[
\lim _{y \to 0+} yg(x^1+\ii y,w) =0
\]
for all $w\in\Pi$, and in particular when $w=x^2$.  Thus $g$ satisfies equation (\ref{limyg}).
\end{proof}
\begin{remark}   \rm
The proof shows that for any $g\in\PP$, if $x^1\in\R$ and 
$\lim_{y\to 0+} yg(x^1+\ii y, z^2)=0$  holds for some $z^2\in\Pi$ then it holds for all $z^2 \in\Pi$. 
\end{remark}
\begin{remark} \rm
The parametrization formula (\ref{paramh}), where $g$ is allowed to be a {\em free} function in the Pick class $\PP$, parametrizes the solutions of a relaxed one-point interpolation problem.  That is, in Theorem \ref{paramP2} one may omit the condition (\ref{limyg}) on $g$ and replace the condition $\nabla h(x) = (\beta, 0)$ by:  $\nabla h(x) = (\beta', 0)$ for some $\beta' \leq \beta$.
\end{remark}
We can now invoke Cayley transformation to obtain a parametrization of functions in $\SS$ with prescribed value and gradient at a facial $B$-point.
\begin{theorem} \label{paramS2}
Let $\tau\in\T\times\D,\ \omega\in\T$ and $\al > 0$.  Suppose that $\tau^1 \neq 1, \ \omega\neq 1$.  The functions $\ph\in\SS$ such that $\tau$ is a $B$-point for $\ph$, that $\ph(\tau) = \omega$ and that $\nabla\ph(\tau)=(\omega\overline{\tau^1}\al, 0)^T$ are precisely the functions of the form
\beq\label{paramph}
\ph(\la) =\frac{h(z) -\ii}{h(z)+\ii}
\eeq
where
\[
z^j = \ii\frac{1+\la^j}{1-\la^j} \quad \mbox{ for } j=1,2,
\]
and
\beq \label{formh}
h(z) = \ii \frac{1+\omega}{1-\omega} + \frac{1}{\frac{1}{\beta(z^1-x^1)}-g(z)}
\eeq
where 
\[ 
   \beta= \frac{1-\re \tau^1}{1-\re \omega} \ \al , \qquad    x^1= \ii \frac{1+\tau^1}{1-\tau^1}
\]
and  $g$ is a function in $\PP$ such that
\beq\label{gcond}
\lim_{y\to 0+} yg\left( \ii\frac{1+\tau^1}{1-\tau^1} +\ii y, \ii \frac{1+\tau^2}{1-\tau^2} \right) = 0.
\eeq
\end{theorem}
\begin{proof}
Consider a function $\ph$ of the form (\ref{paramph}) with $h$ as described.  By Theorem \ref{paramP2} we have $h\in\PP$,
\[
x=(x^1,x^2) = \left(  \ii\frac{1+\tau^1}{1-\tau^1}, \ii \frac{1+\tau^2}{1-\tau^2} \right)
\]
is a $B$-point for $h$,  
\[
h(x) = \ii \frac{1+\omega}{1-\omega} \qquad  \mbox{ and } \qquad\nabla h(x) = \begin{pmatrix} \beta \\ 0 \end{pmatrix}.
\]
It follows that $\tau$ is a $B$-point for $\ph$, $\ph(\tau)=\omega$ and, by equation (\ref{gotnablah}) 
\begin{align*}
\nabla\ph(\tau) &= \begin{pmatrix} \eta^1\\\eta^2\end{pmatrix} =\left(\frac{1-\omega}{1-\tau^1}\right)^2 \nabla h(x)  = \left(\frac{1-\omega}{1-\tau^1}\right)^2 \begin{pmatrix} \beta \\ 0 \end{pmatrix}\\
  &=  \left(\frac{1-\omega}{1-\tau^1}\right)^2 \frac{1-\re\tau^1}{1-\re\omega} \begin{pmatrix} \al \\ 0 \end{pmatrix}  \\
  &= \omega\overline{\tau^1} \begin{pmatrix} \al \\ 0 \end{pmatrix}
\end{align*}
as required.

Conversely, if $\tau$ is a $B$-point for $\ph\in\SS$, if $\ph(\tau) = \omega$ and  $\nabla\ph(\tau)=(\omega\overline{\tau^1}\al, 0)^T$ then $h$ defined by equation (\ref{paramph}) belongs to $\PP$, $x$ is a $B$-point for $h$ and $h(x)=\ii(1+\omega)/(1-\omega), \ \nabla h(x) = (\beta, 0)^T$.  We may therefore apply Theorem \ref{paramP2} to obtain the parametric expression (\ref{formh}) for $h$.
\end{proof}
\begin{remark} {\rm
Again there is a ``relaxed" version of the parametrization.  In Theorem \ref{paramS2}, if one enlarges the class of $\ph$ by replacing the condition  $\nabla\ph(\tau) =(\omega\overline{\tau^1}\al,0)^T$     by $\nabla\ph(\tau) =(\omega\overline{\tau^1}\al',0)^T$ for some $\al'$ with $0 < \al'\leq\al$, then one obtains the same parametrization but without the limit condition (\ref{gcond}) on $g$.}
\end{remark}

%% file: applic.tex
\section{An interpolation problem} \label{applic}
 Theorem \ref{Bfaces} suggests a natural interpolation problem:
to describe the functions $\ph\in\SS$ having $B$-points at finitely many given points  in $\T\times\D$ and $\D\times\T$ and with prescribed values of $\ph$ and $\nabla\ph$ at those points.  In view of Corollary \ref{2faces}, if there are interpolation nodes on both types of face then the target values of $\ph$ must all be the same, but the target values of $\nabla\ph$ may differ.  The parametrization theorems allow us to solve this problem.  We state the result for the half-plane version.   Taking a slight liberty with terminology of D. Sarason, we define\\
{\bf Problem $\partial$NP$\Pick_2$ (facial)}:  {\em Given $x_1,\dots,x_m \in\R\times\Pi, \ y_1, \dots ,y_n\in\Pi\times\R$ for some $m,n \geq 0$ and 
 $\xi\in \R, \ \beta_1,\dots,\beta_m>0, \gamma_1,\dots,\ga_n>0$, determine whether there exists a function $h$ in the two-variable Pick class $\Pick_2$ that satisfies, for $j=1,\dots,m, \ k=1,\dots, n$,
\begin{enumerate}
\item[{\rm(1)}] $x_j$ and $y_k$ are $B$-points for $h$ 
\item[{\rm (2)}] $h(x_j)=\xi=h(y_k)$ and
\item[{\rm (3)}]  \qquad \qquad $
 \nabla h(x_j)=\begin{pmatrix} \beta_j \\ 0 \end{pmatrix}, \qquad \nabla h(y_k) = \begin{pmatrix} 0 \\  \ga_k \end{pmatrix}$.
\end{enumerate}
Describe the set of all such functions $h$ when they exist. }

 We shall write $e_1=(1,0), \ e_2=(0,1)$ for the standard basis in $\C^2$.
\begin{theorem}\label{2interp}
Problem $\partial$NP$\Pick_2$ (facial) always has infinitely many solutions.  The solution set consists of all functions $h$ expressible in the form
\beq\label{allh}
h(z) = \xi + \frac{1} { r(z) -f(z)} \quad \mbox{ for all } z\in\Pi
\eeq
 for some $f\in\Pick_2$ such that 
\beq\label{condf}
\lim_{t\to 0+} tf(x_j+\ii te_1) = 0 = \lim_{t\to 0+} tf(y_k+\ii t e_2), \quad 1\leq j\leq m, \ 1\leq k\leq n,
\eeq
where
\[
r(z) = \sum_{j=1}^m \frac{1}{\beta_j(z^1-x_j^1)} + \sum_{k=1}^n \frac{1}{\ga_k (z^2-y_k^2)}.
\]
In particular, the function $h = \xi + 1/(r-c)$
is a solution of Problem $\partial$NP$\Pick_2$ (facial) for any $c\in\R$.
\end{theorem}
\begin{proof}
We prove necessity by induction on $m+n$.  The assertion holds trivially if $m=n=0$, where empty sums are as usual defined to be $0$.  Suppose $m > 0$ and that necessity holds for $m+n-1$.
Let $h\in\Pick_2$ satisfy conditions (1)-(3).  By Theorem \ref{paramP2}  there exists $g\in\Pick_2$ such that 
\beq \label{limtg}
\lim_{t\to 0+} tg(x_1+\ii te_1) = 0
\eeq
and
\beq\label{gotg}
h(z) = \xi +\frac{1}{\frac{1}{\beta_1 (z^1-x_1^1)}-g(z)} \qquad \mbox{ for all } z\in \Pi.
\eeq
Let $\tilde g = -1/g$; then $\tilde g\in\Pick_2$, and we claim that $x_2,\dots,x_m, y_1,\dots, y_n$ are $B$-points for $\tilde g$.  Indeed, by Theorem \ref{Bfaces}, all these points are $C$-points for $h$, and in view of the hypotheses (2) and (3) on $h$,
\[
h(x_j+\ii te_1) = \xi + \beta_j\ii t +o(t), \quad
h(y_k+\ii te_2) = \xi + \ga_k \ii t +o(t) \quad \mbox{ as } t\to 0+.
\]
Hence, for $j=2,\dots,m$,
\[
 g(x_j+\ii te_1) = -\frac{1}{ \ii \beta_j t + o(t)} + \frac{1}{\beta_1(x_j^1+\ii t -x_1^1)},
\]
from which we easily calculate
\[
\frac{\tilde g(x_j+\ii te_1)}{t} = \beta_j + o(1).
\]
Thus $x_2,\dots,x_m$ are $B$-points for $\tilde g$, and $\tilde g(x_j) =0, \ j= 2\dots,m$.  Likewise $y_1,\dots,y_m$ are $B$-points for $\tilde g$, and $\tilde g$ takes the value $0$ at the $y_k$.  Moreover, by Theorem \ref{Bfaces},
\[
\nabla\tilde g(x_j) = \begin{pmatrix} \beta_j \\0 \end{pmatrix}, \quad
  \nabla\tilde g(y_k) = \begin{pmatrix} 0 \\ \ga_k  \end{pmatrix}
\]
for $j=2,\dots,m, \ k=1,\dots,n$.

By the inductive hypothesis there exists $f\in\Pick_2$ such that
\[
\lim_{t\to 0+} f(x_j + \ii te_1) = 0= \lim_{t\to 0+}f(y_k+\ii te_2)
\]
for $j=2,\dots,m, \ k=1,\dots,n$ and 
\[
\tilde g(z) = 0+ \frac{1}{r_1(z)-f(z)} \qquad \mbox{ for all } z\in\Pi
\]
where
\[
r_1(z)= \sum_{j=2}^m \frac{1}{\beta_j(z^1-x_j^1)} + \sum_{k=1}^n\frac{1}{\ga_k(z^2-y_k^2)}.
\]
Thus
\[
g(z) = -\frac{1}{\tilde g(z)} = f(z) - r_1(z).
\]
This relation in conjunction with equation (\ref{limtg}) tells us that
\[
\lim_{t\to 0+} tf(x_1+\ii te_1) =0,
\]
so that $f$ satisfies conditions (\ref{condf}).
On substituting for $g$ in equation (\ref{gotg}) we obtain 
\begin{align*}
h(z) &= \xi+ \frac{1}{\frac{1}{\beta_1(z^1-x_1^1)} +r_1(z)-f(z)} \\
   &= \xi+ \frac{1}{r(z)-f(z)},
\end{align*}
which is the desired relation (\ref{allh}).  We have proved necessity in Theorem \ref{2interp}.

Conversely, suppose that $h$ is expressible in the form (\ref{allh}), where $f\in\Pick_2$ satisfies conditions (\ref{condf}).   Observe that, for $t>0$,
\begin{align*}
tr(x_j+\ii te_1) &= \sum_{\ell=1}^m \frac{t}{\beta_\ell(x_j^1+\ii t-x_\ell^1)} + \sum_{k=1}^n \frac{t}{\ga_k(x_j^2-y_k^2)} \\
   &= \frac{1}{\beta_j\ii} + o(1).
\end{align*}
Hence
\begin{align*}
\frac{\im h(x_j +\ii te_1)}{t} &= \frac{1}{t} \im \frac{1}{ r(x_j+\ii te_1) -f(x_j+\ii te_1)} \\
   &= \im \frac{1}{ \frac{1}{\beta_j\ii} + o(1) - tf(x_j+\ii te_1)}.
\end{align*}
In view of the relation (\ref{condf}), it follows that
\[
\lim_{t\to 0+} \frac{\im h(x_j+\ii te_1)}{t} = \beta_j, \quad j=1,\dots, m,
\]
so that $x_j$ is a $B$-point for $h$ and $\nabla h(x_j) = (\beta_j, 0)^T$.  Similarly, $y_k$ is a $B$-point for $h$ and $\nabla h(y_k) = (0, \ga_k)^T$ for $k=1,\dots,n$.  
  Hence $h$ satisfies (1)-(3).  In particular we can take $f$ to be a real constant $c$, yielding the solution $h=\xi + 1/(r-c)$.
\end{proof}

\begin{remark} {\rm
We can also consider a relaxed version of the interpolation problem. \\
{\bf Problem $\partial$NP$\Pick_2'$ (facial)}: {\em  As Problem $\partial$NP$\Pick_2$ (facial), but with condition $(3)$ replaced by
\begin{enumerate}
\item[{\rm (3$'$)}]  \qquad  \qquad $   \nabla h(x_j)=\begin{pmatrix} \beta_j' \\ 0 \end{pmatrix}, \qquad \nabla h(y_k) = \begin{pmatrix} 0 \\  \ga_k' \end{pmatrix}$
\end{enumerate}
for some $\beta_j', \ \ga_k'$ such that  $0 < \beta_j' \leq \beta_j, \ 0< \ga_k' \leq \ga_k$. }

The general solution of this relaxed problem is again given by the formula (\ref{allh}), but now without the limit condition (\ref{condf}) on $f$.}
\end{remark}
Do analogous results hold for the polydisk in dimensions greater than 2?  They may well do, but our methods, depending as they do on the use of models, only give analogous statements for those functions which possess models in the sense of Definition \ref{4.4}  (modified appropriately).   In dimensions higher than 2 such functions constitute a class that is strictly smaller than the Schur class \cite{ag90}.  It is often called the {\em Schur-Agler} class.